\documentclass[12pt]{article}

\usepackage{amsmath,amssymb}
\usepackage{theorem}
\usepackage{rotating}

\newtheorem{theorem}{Theorem}[section]
\newtheorem{corollary}[theorem]{Corollary}
\newtheorem{proposition}[theorem]{Proposition}
\newtheorem{facts}{Facts}[section]
\newtheorem{fact}{Fact}[section]

{\theorembodyfont{\rmfamily}\newtheorem{definition}[theorem]{Definition}
\newtheorem{remark}[theorem]{Remark}}

\makeatletter
\DeclareRobustCommand{\qed}{%
  \ifmmode
  \else \leavevmode\unskip\penalty9999 \hbox{}\nobreak\hfill
  \fi
  \quad\hbox{\qedsymbol}}
\newcommand{\openbox}{\leavevmode
  \hbox to.77778em{%
  \hfil\vrule
  \vbox to.675em{\hrule width.6em\vfil\hrule}%
  \vrule\hfil}}
\newcommand{\qedsymbol}{\openbox}
\newenvironment{proof}[1][\proofname]{\par
  \normalfont
  \topsep6\p@\@plus6\p@ \trivlist
  \item[\hskip\labelsep\itshape
    #1\@addpunct{.}]\ignorespaces
}{%
  \qed\endtrivlist
}
\newcommand{\proofname}{Proof}
\makeatother

\numberwithin{equation}{section}

\newcommand{\TTT}{\mathcal T}
\newcommand{\OOO}{\mathcal O}
\newcommand{\III}{\mathcal I}
\newcommand{\Sy}{\mathrm{Sym}}
\newcommand{\Al}{\mathrm{Alt}}

\newcommand{\Hom}{\operatorname{Hom}}
\newcommand{\SU}{\operatorname{SU}}
\newcommand{\SL}{\operatorname{SL}}
\newcommand{\GL}{\operatorname{GL}}
\newcommand{\Sym}{\operatorname{S}}
\newcommand{\ZZ}{\mathbb Z}
\newcommand{\FF}{\mathbb F}
\newcommand{\CC}{\mathbb C}
\newcommand{\sfA}{\mathsf A}
\newcommand{\sfB}{\mathsf B}
\newcommand{\sfC}{\mathsf C}
\newcommand{\sfD}{\mathsf D}
\newcommand{\sfE}{\mathsf E}
\newcommand{\sfF}{\mathsf F}
\newcommand{\sfG}{\mathsf G}
\newcommand{\GrpG}{\Gamma}
\newcommand{\Ghat}{\widehat\Gamma}
\newcommand{\Caff}{C_{\mathrm{aff}}}
\newcommand{\Cfin}{C_{\mathrm{fin}}}
\newcommand{\rep}[1]{\mathsf{#1}}
\newcommand{\Rep}[1]{\mathbf{#1}}
\newcommand{\REP}{\operatorname{R}}

\begin{document}

\title{\textbf{Quantum affine {C}artan matrices, {P}oincar\'e}\\[2mm]
\textbf{series of binary polyhedral groups, and}\\[2mm]
\textbf{reflection representations}\\[2mm]}

\author{\Large Ruedi Suter\\[8mm]
{\textit{Institut f\"ur Theoretische Physik,}}\\
{\textit{Eidgen\"ossische Technische Hochschule Z\"urich,}}\\
{\textit{ETH H\"onggerberg, 8093 Z\"urich, Switzerland}}\\[2mm]
e-mail: \texttt{suter@math.ethz.ch}\\
\vspace*{-0.05truein}\\
}

\maketitle

\section{Introduction}
One of the first results in Lie representation theory is that the symmetric
powers $\Sym^n(\CC^2)$ (for $n\in\ZZ_{\geqslant0}$) of the standard
representation $\CC^2$ of $\SU(2)$ are representatives for the list of
equivalence classes of the irreducible complex representations of $\SU(2)$.
Take a finite subgroup $\GrpG\subseteq\SU(2)$ and let $\rep i$ be an irreducible
complex representation of $\GrpG$. It is a basic question to ask what is the
multiplicity of $\rep i$ in the restriction of $\Sym^n(\CC^2)$ to~$\GrpG$. This
question has been addressed and answered by Kostant \cite{Ko} in a beautiful
way. A crucial ingredient in his approach is a Coxeter transformation
$c_{\mathrm{aff}}$ of the affine Weyl group associated to $\GrpG$ via the McKay
correspondence \cite{McK1,McK2,St}. Kostant writes $c_{\mathrm{aff}}$ as a
product of two involutions $r_1$ and $r_2$ where $r_1$ and $r_2$ themselves are
products of commuting simple reflections. This is only possible if the
affine Coxeter-Dynkin diagram has no odd cycle. Hence type $\sfA_{2n}$ must be
omitted in this approach. In a somewhat different context Springer \cite{Sp}
reproved Kostant's results in his paper in the Mathematische Annalen volume
dedicated to Hirzebruch on his sixtieth birthday. See also \cite{Sp2}. Earlier
papers by Gonzalez-Sprinberg and Verdier \cite{GV} as well as by Kn\"orrer
\cite{Kn} also deal with the question stated above (but their main goal was
something else, namely, the McKay correspondence and singularity theory). They
used that $\GrpG$ is an index~$2$ subgroup of a complex reflection group.

As is often the case in well developed and elementary subjects, one can hardly
avoid rediscovering previously known results. In the first few sections I shall
show how one can derive some invariant theoretic results very easily and quickly
by using quantum affine Cartan matrices. In this largely expository part I
tried to avoid too much overlap with other expositions, and I hope that even
the informed reader will find here some new aspects. Sect.~\ref{hifcrg} displays
for each pair $(\GrpG,\rep i)$ where $\GrpG$ is a (binary) tetrahedral,
octahedral, or icosahedral group and $\rep i$ is an irreducible representation
of $\GrpG$ a homomorphism into a finite complex reflection group $G$ such that
the reflection representation of $G$ restricts to the representation $\rep i$
of $\GrpG$.

There were two main stimuli for writing the present note. The first was McKay's
short paper on semi-affine Coxeter-Dynkin diagrams in the issue of the Canadian
Journal of Mathematics dedicated to Coxeter on his ninetieth birthday
\cite{McK3}, and the second was a recent preprint of Kostant's \cite{Ko2} which
came to me as an inspiration to dig out and complete my notes from 2000/2001.

\mbox{}\\\noindent
\textit{Acknowledgements.} I am pleased to thank J\"urg Fr\"ohlich for
inviting me to visit ETH H\"onggerberg and ITP at ETH Z\"urich for its
hospitality and financial support.

\section{McKay's correspondence}
The starting point is the same as in Kostant's paper: we accept the McKay
correspondence. It establishes a one-to-one correspondence between the set of
conjugacy classes of finite subgroups of $\SU(2)$ and the set of simply-laced
(types $\sfA$, $\sfD$, $\sfE$) affine Coxeter-Dynkin diagrams, and it gives the
character table of such a finite subgroup $\GrpG$ in terms of eigenvectors of
the corresponding affine Cartan matrix. What we actually need from the McKay
correspondence other than the one-to-one correspondence alluded to above is the
following corollary.

\begin{corollary}[Corollary to the McKay correspondence]\ 
Let $\GrpG$ be a finite sub\-group of \,$\SU(2)$. Its unitary dual $\Ghat$ can
then be identified with the set of vertices in the affine Coxeter-Dynkin diagram
associated to $\GrpG$ via the McKay correspondence. (The affine vertex thereby
corresponds to the trivial $1$-dimensional representation.) Let
$\rep{st}=\CC^2|_\GrpG$ (which we will sometimes simply write as $\CC^2$ for
easier readability) denote the standard representation of $\GrpG$, which comes
from the inclusion $\GrpG\subseteq\SU(2)$. For $\rep i\in\Ghat$ the tensor
product $\rep{st}\otimes\rep i$ decomposes as
\begin{equation}\label{mckay}
\rep{st}\otimes\rep i\cong\bigoplus_{\rep j\in N(\rep i)}\rep j
\end{equation}
where $N(\rep i)$ is the (multi-)set of neighbours of $\rep i$, i.\,e., consists
of those vertices in $\Ghat$ that are connected to $\rep i$ by an edge. For
instance, if we write $\rep1$ for the trivial $1$-dimensional representation,
then $\rep{st}\cong\bigoplus\limits_{\rep j\in N(\rep1)}\rep j$. Actually,
$\rep{st}$ is irreducible unless $\GrpG$ is a cyclic group. (And if $\GrpG$
has order $1$ or $2$, then the neighbours count with multiplicity $2$.)
\end{corollary}

\begin{definition}\label{Poincareseries}
For $\rep i\in\Ghat$ we define the generating function (Poincar\'e series)
$$P_{\rep i}(t):=\sum_{n=0}^\infty\dim\Hom_\GrpG\bigl(\rep i,\Sym^n(\CC^2)\bigr)
\cdot t^n$$
whose $t^n$ coefficient is the multiplicity of $\rep i$ in $\Sym^n(\CC^2)$.
\end{definition}

Consider the classical $\SU(2)$ (or its restriction to $\GrpG$) Clebsch-Gordan
decomposition
\begin{equation}\label{ClebschGordan}
\CC^2\otimes\Sym^n(\CC^2)\cong\Sym^{n+1}(\CC^2)\oplus\Sym^{n-1}(\CC^2)
\end{equation}
for $n\in\ZZ_{\geqslant0}$. In an explicit realization as a space of
polynomials, we can think of $\CC^2$ as a space of linear polynomials in two
variables and of $\Sym^n(\CC^2)$ as a space of homogeneous binary polynomials of
degree $n$. So on the left side we see a space of bihomogeneous polynomials in
$2+2$ variables of bidegree $(1,n)$ (and in particular of total degree $n+1$).
So also the space on the right side can be realized by the same space of
bihomogeneous polynomials. Taking into account the degrees and using the McKay
correspondence, we arrive at the equation
\begin{equation}\label{eq1}
(1+t^2)P_{\rep i}(t)-\delta_{\rep i}=t\sum_{\rep j\in N(\rep i)}P_{\rep j}(t)
\end{equation}
with $\delta_{\rep1}=1$ and $\delta_{\rep i}=0$ for $\rep1\neq\rep i\in\Ghat$.
If you prefer to see a step-by-step derivation of formula~(\ref{eq1}), here is
one:
\begin{align}\notag
\lefteqn{(1+t^2)P_{\rep i}(t)-\delta_{\rep i}}\qquad&\\\notag
&=\sum_{n=0}^\infty\dim\Hom_\GrpG\bigl(\rep i,\Sym^n(\CC^2)\bigr)\cdot t^n
-\delta_{\rep i}+\sum_{n=0}^\infty\dim\Hom_\GrpG\bigl(\rep i,\Sym^n(\CC^2)\bigr)
\cdot t^{n+2}\\\notag
&=t\left(\sum_{n=0}^\infty\dim\Hom_\GrpG\bigl(\rep i,\Sym^{n+1}(\CC^2)\bigr)
\cdot t^n+\sum_{n=0}^\infty\dim\Hom_\GrpG\bigl(\rep i,\Sym^{n-1}(\CC^2)\bigr)
\cdot t^n\right)\\\notag
&=t\sum_{n=0}^\infty\dim\Hom_\GrpG\bigl(\rep i,\Sym^{n+1}(\CC^2)
\oplus\Sym^{n-1}(\CC^2)\bigr)\cdot t^n\\\notag
&=t\sum_{n=0}^\infty\dim\Hom_\GrpG\bigl(\rep i,\CC^2\otimes\Sym^n(\CC^2)\bigr)
\cdot t^n
\intertext{and since $\rep{st}$ is isomorphic to its dual representation we can
continue}\notag
&=t\sum_{n=0}^\infty\dim\Hom_\GrpG\bigl(\rep{st}\otimes\rep i,\Sym^n(\CC^2)
\bigr)\cdot t^n=t\sum_{n=0}^\infty\dim\Hom_\GrpG\Bigl({\textstyle
\bigoplus\limits_{\rep j\in N(\rep i)}}\rep j\,,\Sym^n(\CC^2)\Bigr)
\cdot t^n\\\tag{\ref{eq1}}
&=t\sum_{n=0}^\infty\sum_{\rep j\in N(\rep i)}\dim\Hom_\GrpG\bigl(
\rep j,\Sym^n(\CC^2)\bigr)\cdot t^n=t\sum_{\rep j\in N(\rep i)}P_{\rep j}(t).
\end{align}
The equation (\ref{eq1}) when written as a matrix equation looks like
\begin{equation}\label{eq2}
\Caff(t)P(t)=\delta.
\end{equation}
Here $\Caff(t)$ is the \emph{quantum affine Cartan matrix} for $\GrpG$. Its row
and column indices run over the vertices of the affine Coxeter-Dynkin diagram as
for the usual affine Cartan matrix $\Caff$. Now $\Caff(t)$ is gotten from
$\Caff$ by replacing the diagonal elements $2$ by $1+t^2$ and by multiplying the
off-diagonal entries by $t$. [For type $\sfA_0$ we have
$\Caff(t)=\bigl((1-t)^2\bigr)$.] The vectors $P(t)$ and $\delta$ have
$P_{\rep i}(t)$ and $\delta_{\rep i}$ as entries, respectively. Note that
$\Caff(t)\equiv\mathbf 1\pmod t$, so that the matrix $\Caff(t)$ is invertible in
the power series ring, and the entries of $\Caff(t)^{-1}$ can obviously be
written as rational functions in $t$. So we obtain the various Poincar\'e series
$P_{\rep i}(t)$ simply by solving (\ref{eq2}) for $P(t)$, that is,
$P(t)=\Caff(t)^{-1}\delta$.

The following theorem summarizes what we have gotten so far.

\begin{theorem}\label{ThmPS}
Let $\GrpG$ be a finite subgroup of \,$\SU(2)$. We identify its unitary dual
$\Ghat$ with the vertices of an affine Coxeter-Dynkin diagram with affine
Cartan matrix $\Caff$. Then the vector of Poincar\'e series
$P(t)=\bigl(P_{\rep i}(t)\bigr)_{\rep i\in\Ghat}$ in
Definition~\textup{\ref{Poincareseries}} is
$$P(t)=\Caff(t)^{-1}\delta$$
where $\Caff(t)=(1-t)^2\mathbf1+t\Caff$ is the quantum affine Cartan matrix and
$\delta=(\delta_{\rep i})_{\rep i\in\Ghat}$ is the vector such that
$\delta_{\rep i}=1$ if $\rep i$ is the trivial representation and
$\delta_{\rep i}=0$ otherwise.
\end{theorem}

Cramer's rule yields the following corollary.

\begin{corollary}\label{psinv}
The Poincar\'e series for the invariant ring $\Sym^\ast(\CC^2)^\GrpG$ is
\begin{equation}\label{psinveq}
P_{\rep1}(t)=\frac{\det\Cfin(t)}{\det\Caff(t)}\ .
\end{equation}
Here $\Cfin(t)$ is the submatrix of $\Caff(t)$ with row and column indices
different from $\rep1$.
\end{corollary}

The most natural way to derive the crucial equation (\ref{eq1}) that led us to
introduce the quantum affine Cartan matrix actually comes from the following
general formula that takes place in the formal power series ring with
coefficients in the complex representation ring $\REP(\GrpG)$ of our finite
group $\GrpG$. (By some slight abuse of notation we denote representations or
their equivalence classes in the representation ring by the same letters.)
\begin{equation}\label{lambdasigma}
\frac{\rep1}{\lambda_{-t}(\rep i)}=\sigma_t(\rep i)\in\REP(\GrpG)[\![t]\!]
\end{equation}
Here $\lambda_{-t}(\rep i)$ and $\sigma_t(\rep i)$ are defined as usual, namely,
if $\rep i\in\REP(\GrpG)$ is the class of any (honest) finite-dimensional
complex representation of $\GrpG$, then we look at its exterior powers
$\bigwedge{}^{\!k}\,\rep i$ and symmetric powers $\Sym^n(\rep i)$ and put
\begin{align*}
\lambda_{-t}(\rep i)&=\sum_{k=0}^{\dim\rep i}{\textstyle\bigwedge{}^{\!k}\,
\rep i}\cdot t^k\quad\mbox{and}\quad
\sigma_t(\rep i)=\sum_{n=0}^\infty\Sym^n(\rep i)\cdot t^n.
\end{align*}
For $\rep i=\rep{st}$ we have
$\lambda_{-t}(\rep{st})=\rep1-\rep{st}\,t+\rep1 t^2$ and since
$\sigma_t(\rep{st})=\sum\limits_{\rep i\in\Ghat}\rep i\cdot P_{\rep i}(t)$ the
identity (\ref{lambdasigma}) says that
$$(\rep1-\rep{st}\,t+\rep1 t^2)\Bigl(\sum_{\rep i\in\Ghat}\rep i\cdot
P_{\rep i}(t)\Bigr)=\rep1.$$
Inserting McKay's formula (\ref{mckay}) we get
$$(1+t^2)\Bigl(\sum_{\rep i\in\Ghat}\rep i\cdot P_{\rep i}(t)\Bigr)-\rep1=
t\sum_{\rep j\in\Ghat}\sum_{\rep i\in N(\rep j)}\rep i\cdot P_{\rep j}(t).$$
Its $\rep i$th component is just equation (\ref{eq1}) because
$\rep i\in N(\rep j)$ is equivalent to $\rep j\in N(\rep i)$. An obvious
advantage of formula~(\ref{lambdasigma}) in comparison with the previous
\textit{ad hoc} approach to derive formula (\ref{eq1}) is that we can use
formula (\ref{lambdasigma}) for any representation and not just for the
standard $2$-dimensional representation. We will use it in Sect.~\ref{sect}.

\section{An example---$\sfE_8$}
In this short section we shall apply the formula in Corollary~\ref{psinv} to
compute the Poincar\'e series of the invariant ring of the binary icosahedral
group in its standard $2$-dimensional representation. (Another way of computing
the Poincar\'e series would be to use the classical Molien formula.) The binary
icosahedral group corresponds to $\sfE_8$ via the McKay correspondence, and its
quantum affine Cartan matrix reads

\renewcommand{\arraystretch}{1.4}
{\footnotesize
$$\begin{pmatrix}
1+t^2&-t\phantom{-}&0&0&0&0&0&0&0\\
-t\phantom{-}&1+t^2&-t\phantom{-}&0&0&0&0&0&0\\
0&-t\phantom{-}&1+t^2&-t\phantom{-}&0&0&0&0&0\\
0&0&-t\phantom{-}&1+t^2&-t\phantom{-}&0&0&0&0\\
0&0&0&-t\phantom{-}&1+t^2&-t\phantom{-}&0&0&0\\
0&0&0&0&-t\phantom{-}&1+t^2&-t\phantom{-}&-t\phantom{-}&0\\
0&0&0&0&0&-t\phantom{-}&1+t^2&0&-t\phantom{-}\\
0&0&0&0&0&-t\phantom{-}&0&1+t^2&0\\
0&0&0&0&0&0&-t\phantom{-}&0&1+t^2
\end{pmatrix}
\quad
\setlength{\unitlength}{7mm}
{
\begin{picture}(4,0)(-2,3.3)
\multiput(1,0)(0,1){8}{\circle{0.2}}
\multiput(1,0)(0,1){8}{\circle{0.18}}
\put(0,2){\circle{0.2}}
\put(0,2){\circle{0.18}}
\multiput(1,0.1)(0,1){7}{\line(0,1){0.8}}
\put(0.1,2){\line(1,0){0.8}}
\put(1.3,7){\makebox(0,0)[l]{$1$}}
\put(1.3,6){\makebox(0,0)[l]{$2$}}
\put(1.3,5){\makebox(0,0)[l]{$3$}}
\put(1.3,4){\makebox(0,0)[l]{$4$}}
\put(1.3,3){\makebox(0,0)[l]{$5$}}
\put(1.3,2){\makebox(0,0)[l]{$6$}}
\put(1.3,1){\makebox(0,0)[l]{$7$}}
\put(1.3,0){\makebox(0,0)[l]{$9$}}
\put(-0.3,2){\makebox(0,0)[r]{$8$}}
\end{picture}
}
$$
}%
\renewcommand{\arraystretch}{1}%
where the rows and columns are ordered according to the numbering of the
vertices in the affine $\sfE_8$ diagram displayed above.

One computes
\begin{align*}
\det\Cfin(t)&=1+t^2-t^6-t^8-t^{10}+t^{14}+t^{16},\\
\det\Caff(t)&=1+t^2-t^6-t^8-t^{10}-t^{12}+t^{16}+t^{18}.
\end{align*}
The quotient can be put into the form
$$P_{\rep1}(t)=\frac{\det\Cfin(t)}{\det\Caff(t)}
=\frac{1+t^{30}}{(1-t^{12})(1-t^{20})}\ .$$
This a a classical result known from the theory of Kleinian singularities.

\begin{fact}
Let $\GrpG$ be a finite subgroup of \,$\SU(2)$. The Poincar\'e series for the
invariant ring $\Sym^\ast(\CC^2)^\GrpG$ is
\begin{equation}\label{eqabh}
P_{\rep1}(t)=\frac{1+t^h}{(1-t^a)(1-t^b)}
\end{equation}
where $h=\sum\limits_{\rep i\in\Ghat}\dim\rep i$ is the Coxeter number,
$a=2\max\bigl\{\dim\rep i\bigm|\rep i\in\Ghat\bigr\}$, and $b=h+2-a$, and it
turns out that $ab=2|\GrpG|$.
\end{fact}

In general, Kostant \cite{Ko} proved that one can write the rational functions
$P_{\rep i}(t)$, which we expressed as matrix entries of the inverse
quantum affine Cartan matrix in Theorem~\ref{ThmPS}, as
$$P_{\rep i}(t)=\frac{z_{\rep i}(t)}{(1-t^a)(1-t^b)}$$
for certain polynomials $z_{\rep i}(t)$ with nonnegative integer coefficients
and that these polynomials can be described in a beautiful way by considering
the action of a finite Coxeter transformation $c_{\mathrm{fin}}$ on the root
system and by intersecting the orbits with the set of positive roots that are
not perpendicular to the highest root. The details are in \cite{Ko,Ko2}.

To finish the $\sfE_8$ example, let us display the polynomials
$z_{\rep i}(t)$.

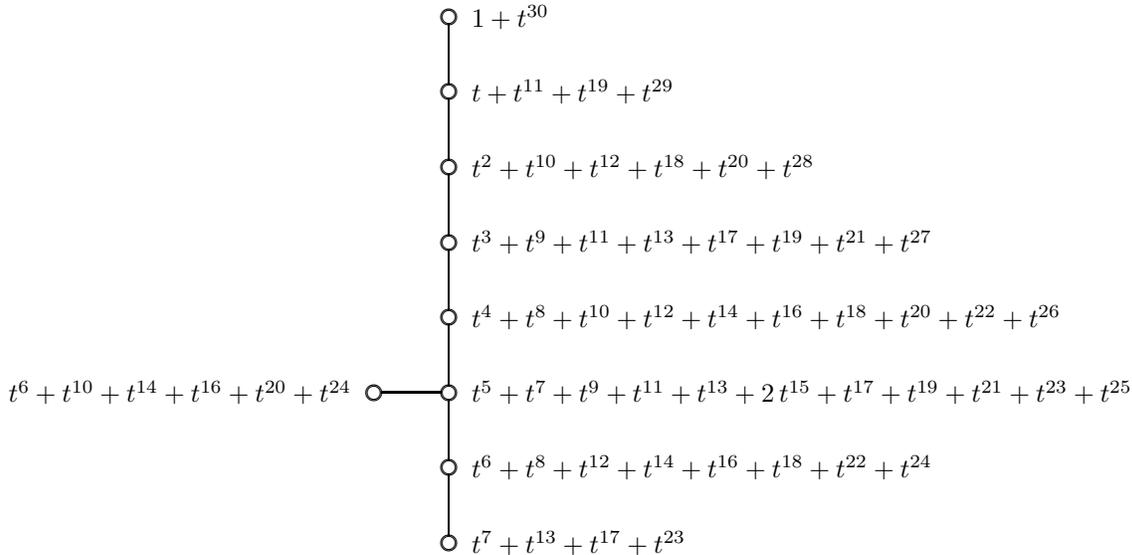
\setlength{\unitlength}{1cm}
\begin{center}
{\footnotesize
\begin{picture}(15,7.5)(-5,0)
\multiput(1,0)(0,1){8}{\circle{0.2}}
\multiput(1,0)(0,1){8}{\circle{0.18}}
\put(0,2){\circle{0.2}}
\put(0,2){\circle{0.18}}
\multiput(1,0.1)(0,1){7}{\line(0,1){0.8}}
\put(0.1,2){\line(1,0){0.8}}
\put(1.3,7){\makebox(0,0)[l]{$1+t^{30}$}}
\put(1.3,6){\makebox(0,0)[l]{$t+t^{11}+t^{19}+t^{29}$}}
\put(1.3,5){\makebox(0,0)[l]{$t^2+t^{10}+t^{12}+t^{18}+t^{20}+t^{28}$}}
\put(1.3,4){\makebox(0,0)[l]{$t^3+t^9+t^{11}+t^{13}+t^{17}+t^{19}+
t^{21}+t^{27}$}}
\put(1.3,3){\makebox(0,0)[l]{$t^4+t^8+t^{10}+t^{12}+t^{14}+t^{16}+t^{18}+
t^{20}+t^{22}+t^{26}$}}
\put(1.3,2){\makebox(0,0)[l]{$t^5+t^7+t^9+t^{11}+t^{13}+2\,t^{15}+t^{17}+
t^{19}+t^{21}+
t^{23}+t^{25}$}}
\put(1.3,1){\makebox(0,0)[l]{$t^6+t^8+t^{12}+t^{14}+t^{16}+t^{18}+t^{22}+
t^{24}$}}
\put(1.3,0){\makebox(0,0)[l]{$t^7+t^{13}+t^{17}+t^{23}$}}
\put(-0.3,2){\makebox(0,0)[r]{$t^6+t^{10}+t^{14}+t^{16}+t^{20}+t^{24}$}}
\end{picture}
}
\end{center}

\section{A numerological table}
Up to now we have looked at the parameters $a,b,h$ for $\sfA\sfD\sfE$ types.
Some of the corresponding affine Coxeter-Dynkin diagrams can be folded into
affine Coxeter-Dynkin diagrams of $\sfC\sfB\sfF\sfG$ types. In this way
$\sfC_l$ unfolds to $\sfA_{2l-1}$, $\sfB_l$ unfolds to $\sfD_{l+1}$, $\sfF_4$
unfolds to $\sfE_6$, and $\sfG_2$ unfolds to $\sfD_4$. We define the parameters
$a,b,h$ for the folded types to be equal to those parameters for the
corresponding unfolded types. The formula
\begin{equation}\label{P1}
\frac{\det\Cfin(t)}{\det\Caff(t)}=\frac{1+t^h}{(1-t^a)(1-t^b)}
\end{equation}
which for the $\sfA\sfD\sfE$ types follows
from (\ref{psinveq}) and (\ref{eqabh}) is actually valid for all types.
In fact, from a computational point of view
this may be obvious by taking into account symmetries.
The parameters $p,q,r$ (ordered appropriately) in the table are determined
by Proposition~\ref{DetqaC}.

\newlength{\abpqr}\settowidth{\abpqr}{$\frac12(l+1)$}
$$\begin{array}{|l|c|c|c|c|c|c|}\hline
\multicolumn{1}{|c|}{\mbox{type}}&a&b&h&p&q&r\\\hline\hline
\rlap{$\sfA_l$}\qquad(l\geqslant0)&\makebox[\abpqr][c]{$2$}&
\makebox[\abpqr][c]{$l+1$}&\makebox[\abpqr][c]{$l+1$}&\frac12(l+1)&\frac12(l+1)
&\makebox[\abpqr][c]{$1$}\\
\rlap{$\sfD_l$}\qquad(l\geqslant4)&4&2l-4&2l-2&l-2&2&2\\
\sfE_6&6&8&12&3&3&2\\
\sfE_7&8&12&18&4&3&2\\
\sfE_8&12&20&30&5&3&2\\\hline
\rlap{$\sfC_l$}\qquad(l\geqslant2)&2&2l&2l&l&1&1\\
\rlap{$\sfB_l$}\qquad(l\geqslant3)&4&2l-2&2l&l-1&2&1\\
\sfF_4&6&8&12&3&2&1\\
\sfG_2&4&4&6&2&1&1\\\hline
\end{array}$$

\begin{proposition}\label{DetqaC}
The determinant of a quantum affine Cartan matrix $\Caff(t)$
is
\begin{equation}\label{DetqaCeq}
\det\Caff(t)=\frac{(1-t^{2p})(1-t^{2q})(1-t^{2r})}{1-t^2}
\end{equation}
where $p,q,r$ are given in the table.
\end{proposition}
\begin{proof}
The proof is an elementary computation. (See also \cite{LT}.)
\end{proof}

Looking at the degree, we obtain the following corollary.
\begin{corollary}
$p+q+r=l+2$.
\end{corollary}

The numbers $p$, $q$, and $r$ have the usual interpretation for the
$\sfA\sfD\sfE$ types. For instance, $p-1$, $q-1$, and $r-1$ are the arm lengths
in the $\mathsf Y$-shaped \emph{finite} Coxeter-Dynkin diagram (to be
interpreted with a grain of salt in case of type $\sfA_{2n}$).

The next corollary gives a formula for the determinant of the finite Cartan
matrix $\Cfin=\Cfin(1)$. It generalizes Saito's relation \cite{Sa1} to all
types.
\begin{corollary}
$\displaystyle\det\Cfin=\frac{8\,pqr}{ab}$\ .
\end{corollary}
\begin{proof}
Solve (\ref{P1}) and (\ref{DetqaCeq}) for $\Cfin(t)$,
cancel the factors $(1-t)^3$ in the numerator and denominator, and evaluate at
$t=1$.
\end{proof}

\begin{remark}
Using the fact that $ab=2|\GrpG|$ and the formula for $|\GrpG|$ in terms of
$p,q,r$ (see the last entry of Facts~\ref{facts} below), one writes $\det\Cfin$
in terms of $p,q,r$ for the $\sfA\sfD\sfE$ types. It turns out that the same
expression in $p,q,r$ is a well-known quantity for the $\sfC\sfB\sfF\sfG$ types,
too.
$$pq+qr+pr-pqr=\begin{cases}
\det\Cfin&\mbox{for the $\sfA\sfD\sfE$ types},\\
l+1&\mbox{for the $\sfC\sfB\sfF\sfG$ types}.
\end{cases}$$
\end{remark}

Before stating the next corollary let me recall how the spectrum of the Cartan
matrix $\Cfin$ or $\Caff$ is related to the spectrum of a corresponding Coxeter
transformation $c_{\mathrm{fin}}$ or $c_{\mathrm{aff}}$ (see \cite{BLM,Co}).
Namely, there are complex numbers $\chi_1,\dots,\chi_n$ satisfying
\mbox{$\chi_j+\chi_{n+1-j}=2\pi i$} (for $j=1,\dots,n$) such that
$e^{\chi_1},\dots,e^{\chi_n}$ are the eigenvalues of a Coxeter transformation.
The eigenvalues of the corresponding Cartan matrix are then
$2-e^{\chi_1/2}-e^{-\chi_1/2},\dots,2-e^{\chi_n/2}-e^{-\chi_n/2}$.
This is true if the underlying Coxeter-Dynkin diagram is a tree (or forest),
which covers all finite cases and all irreducible affine cases except affine
$\sfA$ type. In the latter case the Coxeter-Dynkin diagram is a cycle with $l+1$
vertices, and the $(l+1)!$ possible Coxeter transformations occur in
$\lfloor(l+1)/2\rfloor$ spectral classes (for $l\geqslant1$). If
$s_1,\dots,s_{l+1}$ are the simple reflections (still for affine $\sfA_l$ type),
then one can compute (see \cite{Co}) that
$$\prod_{\sigma\in\Sy_{l+1}}\det\bigl(T\cdot\operatorname{id}-s_{\sigma(1)}
\cdots s_{\sigma(l+1)}\bigr)=\prod_{k=1}^l(1-T^k)^{2(l+1)A(l,k)}$$
where $\Sy_{l+1}$ is the symmetric group of degree $l+1$ and
$A(l,k)=\sum\limits_{j=0}^k(-1)^j\,\binom{l+1}{j}\,(k-j)^l$ are
Eulerian numbers.

Coming back to the tree case we evaluate the characteristic polynomial
(in $t^2$) of a Coxeter transformation as
\begin{align*}
\prod_{j=1}^n\bigl(t^2-e^{\chi_j}\bigr)&
=\prod_{j=1}^n\bigl((t-e^{\chi_j/2})(t+e^{\chi_j/2})\bigr)
=\prod_{j=1}^n(t-e^{\chi_j/2})\prod_{j=1}^n(t+e^{\chi_{n+1-j}/2})\\
&=\prod_{j=1}^n\bigl((t-e^{\chi_j/2})(t-e^{-\chi_j/2})\bigr)
\tag{\mbox{using the symmetry $\chi_j+\chi_{n+1-j}=2\pi i$}}\\
&=\prod_{j=1}^n\bigl(t^2-te^{\chi_j/2}-te^{-\chi_j/2}+1\bigr)
=\prod_{j=1}^n\bigl((1-t)^2+t(2-e^{\chi_j/2}-e^{-\chi_j/2})\bigr)
\end{align*}
so that we get
$\det\Cfin(t)=\det(t^2\cdot\operatorname{id}-c_{\mathrm{fin}})$ and
$\det\Caff(t)=\det(t^2\cdot\operatorname{id}-c_{\mathrm{aff}})$ except for
affine $\sfA$ type.

The next corollary gives a characterization of the exponents $m_j$ of a finite
Weyl group. (Recall that $1\leqslant m_j<h$ and $\exp(2\pi i m_j/h)$
(for $j=1,\dots,l$) are the eigenvalues of a finite Coxeter transformation
$c_{\mathrm{fin}}$.)
\begin{corollary}
Setting $q=e^{2\pi i/h}$, we have
$$\prod_{j=1}^l\bigl(t^2-q^{m_j}\bigr)=\frac{(1-t^{2h})(1-t^{2p})(1-t^{2q})
(1-t^{2r})}{(1-t^2)(1-t^a)(1-t^b)(1-t^h)}\ .$$
\end{corollary}
\begin{proof}
The left hand side is $\Cfin(t)$ by the discussion above, and by (\ref{P1})
and (\ref{DetqaCeq}) $\Cfin(t)$ can be expressed by the right hand side.
\end{proof}


\section{Some facts about binary polyhedral groups}
\begin{facts}\label{facts}
Let $\GrpG=\bigl\langle\alpha,\beta,\gamma\bigm|\alpha^p=\beta^q=\gamma^r
=\alpha\beta\gamma\bigr\rangle$ be a (finite) binary polyhedral group with
$p\geqslant q\geqslant r\geqslant1$ and $p=q$ if $r=1$.
\begin{itemize}
\item The collection of parameters $(p,q,r)$ as above is
$$\begin{array}{cl}
(p,p,1)&\mbox{cyclic group of order $2p$,}\\
(p,2,2)&\mbox{binary dihedral (\,= generalized quaternion = dicyclic) group of
order $4p$,}\\
(3,3,2)&\mbox{binary tetrahedral group $\TTT$ (of order $24$),}\\
(4,3,2)&\mbox{binary octahedral group $\OOO$ (of order $48$),}\\
(5,3,2)&\mbox{binary icosahedral group $\III$ (of order $120$).}
\end{array}$$
\item If $\GrpG$ is not cyclic, then the centre $Z(\GrpG)$ has order $2$ and is
generated by $\alpha\beta\gamma$.
\item The conjugacy classes other than $\{1\}$ and $\{\alpha\beta\gamma\}$ are
the following:
\begin{itemize}
\item[$\circ$] $[\alpha^j]$ of size $\bigl|[\alpha^j]\bigr|=|\GrpG|/2p$ for
$1\leqslant j\leqslant p-1$,
\item[$\circ$] $[\beta^j]$ of size $\bigl|[\beta^j]\bigr|=|\GrpG|/2q$ for
$1\leqslant j\leqslant q-1$,
\item[$\circ$] $[\gamma^j]$ of size $\bigl|[\gamma^j]\bigr|=|\GrpG|/2r$ for
$1\leqslant j\leqslant r-1$.
\end{itemize}
Therefore the order of $\Gamma$ is the following well-known expression.
$$|\Gamma|=\frac{4}{\frac1p+\frac1q+\frac1r-1}$$
\end{itemize}
\end{facts}

\section{Explicit computations: the primitive cases $\sfE_6,\sfE_7,\sfE_8$}
Some vertices in the following Coxeter-Dynkin diagrams are decorated. The cross
is for the trivial representation. The spinorial representations, i.\,e., those
representations where the nontrivial central element
$\overline1=\alpha\beta\gamma$ acts by $-1$, are marked by a dot.

\subsection*{Tetrahedral case
$\TTT=\bigl\langle\alpha,\beta,\gamma\bigm|\alpha^3=\beta^3=\gamma^2
=\alpha\beta\gamma\bigr\rangle\cong\SL_2(\FF_3)$}
\begin{remark}
$\alpha\mapsto\binom{-1\ {-1}}{\phantom{+}0\ {-1}}$,
$\beta\mapsto\binom{\phantom{+}1\ {\phantom{+}1}}{-1\ {\phantom{+}0}}$,
$\gamma\mapsto\binom{\,0\ {-1}}{\,1 \ {\phantom{+}0}}$ gives
an isomorphism $\TTT\cong\SL_2(\FF_3)$.
\end{remark}

\setlength{\unitlength}{1cm}
\begin{center}
{\footnotesize
\begin{picture}(4,3)(-2,-2)
\multiput(-2,0)(1,0){5}{\circle{0.2}}
\multiput(-1,0)(2,0){2}{\circle*{0.08}}
\multiput(0,-2)(0,1){2}{\circle{0.2}}
\put(0,-1){\circle*{0.08}}
\put(0,-2){\makebox(0,0){\tiny$\times$}}
\multiput(-1.9,0)(1,0){4}{\line(1,0){0.8}}
\multiput(0,-1.9)(0,1){2}{\line(0,1){0.8}}
\put(-2,0.237){\makebox(0,0)[b]{$\substack{\textstyle[\alpha]\\[0.456ex]
\textstyle\rep{1'}}$}}
\put(-1,0.237){\makebox(0,0)[b]{$\substack{\textstyle[\alpha^2]\\[0.456ex]
\textstyle\rep{2'}}$}}
\put(0,0.237){\makebox(0,0)[b]{$\substack{\textstyle\{\overline1\}\\[0.456ex]
\textstyle\rep3}$}}
\put(1,0.237){\makebox(0,0)[b]{$\substack{\textstyle[\beta^2]\\[0.456ex]
\textstyle\rep{2''}}$}}
\put(2,0.237){\makebox(0,0)[b]{$\substack{\textstyle[\beta]\\[0.456ex]
\textstyle\rep{1''}}$}}
\put(0.237,-1){\makebox(0,0)[l]{$\rep{2}\ [\gamma]$}}
\put(0.237,-2){\makebox(0,0)[l]{$\rep{1}\ \{1\}$}}
\end{picture}
}\end{center}

\newlength{\cwth}\setlength{\cwth}{8mm}

$\begin{array}[b]{r|ccccccc}
\multicolumn{1}{p{\cwth}}{}&\multicolumn{1}{p{\cwth}}{}&
\multicolumn{1}{p{\cwth}}{}&\multicolumn{1}{p{\cwth}}{}&
\multicolumn{1}{p{\cwth}}{}&\multicolumn{1}{p{\cwth}}{}&
\multicolumn{1}{p{\cwth}}{}&\multicolumn{1}{p{\cwth}}{}\\
\substack{\textrm{character}\\\textrm{table}\hfill}
&\{1\}&[\gamma]
&\{\overline1\}&[\alpha^2]\phantom{^2}
&[\beta^2]\phantom{^2}&[\alpha]\phantom{^2}
&[\beta]\phantom{^2}\\\hline
\rep{1\phantom{''}}&1&1&1&1\phantom{^2}&1\phantom{^2}&1\phantom{^2}
&1\phantom{^2}\\
\rep{2\phantom{''}}&2&0&-2\phantom{-}&-1\phantom{^2}\phantom{-}
&-1\phantom{^2}\phantom{-}&1\phantom{^2}&1\phantom{^2}\\
\rep{3\phantom{''}}&3&-1\phantom{-}&3&0\phantom{^2}&0\phantom{^2}
&0\phantom{^2}&0\phantom{^2}\\
\rep{2{}'\phantom{'}}&2&0&-2\phantom{-}&-\rho\phantom{^2}\phantom{-}
&-\rho^2\phantom{-}&\rho^2&\rho\phantom{^2}\\
\rep{2{}''}&2&0&-2\phantom{-}&-\rho^2\phantom{-}&-\rho\phantom{^2}\phantom{-}
&\rho\phantom{^2}&\rho^2\\
\rep{1{}'\phantom{'}}&1&1&1&\rho\phantom{^2}&\rho^2&\rho^2&\rho\phantom{^2}\\
\rep{1{}''}&1&1&1&\rho^2&\rho\phantom{^2}&\rho\phantom{^2}&\rho^2
\end{array}$\quad\ $\rho=\frac{-1+\sqrt{-3}}{2}=\exp(2\pi i/3)$

\subsection*{Octahedral case
$\OOO=\bigl\langle\alpha,\beta,\gamma\bigm|\alpha^4=\beta^3=\gamma^2
=\alpha\beta\gamma\bigr\rangle$}
\begin{remark}
Note that $\OOO\ncong\SL_2(\FF_4)$ (in fact, $\SL_2(\FF_4)\cong\Al_5$ is the
simple group of order~$60$), and also that $\OOO\ncong\SL_2(\ZZ/4\ZZ)$ (in fact,
the group $\SL_2(\ZZ/4\ZZ)$ (of order $48$) has more than one element of order
$2$). But one can embed $\OOO$ into $\SL_2(\FF_7)$, e.\,g., by
$\alpha\mapsto\binom{\,1\ \ 1}{\,2\ \ 3}$,
$\beta\mapsto\binom{-1\ {-3}}{\phantom{+}1\ {\phantom{+}2}}$,
$\gamma\mapsto\binom{\,0\ {-1}}{\,1 \ {\phantom{+}0}}$.
Note that $\operatorname{tr}\binom{\,1\ \ 1}{\,2\ \ 3}=4$ is a square root of
$2$ in $\FF_7$.
\end{remark}

\setlength{\unitlength}{1cm}
\begin{center}
{\footnotesize
\begin{picture}(6,2)(-3,-1)
\multiput(-3,0)(1,0){7}{\circle{0.2}}
\multiput(-2,0)(2,0){3}{\circle*{0.08}}
\put(0,-1){\circle{0.2}}
\put(3,0){\makebox(0,0){\tiny$\times$}}
\multiput(-2.9,0)(1,0){6}{\line(1,0){0.8}}
\put(0,-0.9){\line(0,1){0.8}}
\put(-3,0.237){\makebox(0,0)[b]{$\substack{\textstyle[\alpha]\\[0.456ex]
\textstyle\rep{1'}}$}}
\put(-2,0.237){\makebox(0,0)[b]{$\substack{\textstyle[\alpha^2]\\[0.456ex]
\textstyle\rep{2'}}$}}
\put(-1,0.237){\makebox(0,0)[b]{$\substack{\textstyle[\alpha^3]\\[0.456ex]
\textstyle\rep{3'}}$}}
\put(0,0.237){\makebox(0,0)[b]{$\substack{\textstyle\{\overline1\}\\[0.456ex]
\textstyle\rep{4}}$}}
\put(1,0.237){\makebox(0,0)[b]{$\substack{\textstyle[\beta^2]\\[0.456ex]
\textstyle\rep{3}}$}}
\put(2,0.237){\makebox(0,0)[b]{$\substack{\textstyle[\beta]\\[0.456ex]
\textstyle\rep{2}}$}}
\put(3,0.237){\makebox(0,0)[b]{$\substack{\textstyle\{1\}\\[0.456ex]
\textstyle\rep{1}}$}}
\put(0.237,-1){\makebox(0,0)[l]{$\rep{2''}\ [\gamma]$}}
\end{picture}
}
\end{center}

$\begin{array}[b]{r|cccccccc}
\multicolumn{1}{p{\cwth}}{}&\multicolumn{1}{p{\cwth}}{}&
\multicolumn{1}{p{\cwth}}{}&\multicolumn{1}{p{\cwth}}{}&
\multicolumn{1}{p{\cwth}}{}&\multicolumn{1}{p{\cwth}}{}&
\multicolumn{1}{p{\cwth}}{}&\multicolumn{1}{p{\cwth}}{}&
\multicolumn{1}{p{\cwth}}{}\\
\substack{\textrm{character}\\\textrm{table}\hfill}
&\{1\}&[\beta]
&[\beta^2]&\{\overline1\}
&[\alpha^3]&[\gamma]
&[\alpha^2]&[\alpha]\\\hline
\rep{1\phantom{''}}&1&1&1&1&1&1&1&1\\
\rep{2\phantom{''}}&2&1&-1\phantom{-}&-2\phantom{-}&-\sigma\phantom{-}&0&0
&\sigma\\
\rep{3\phantom{''}}&3&0&0&3&1&-1\phantom{-}&-1\phantom{-}&1\\
\rep{4\phantom{''}}&4&-1\phantom{-}&1&-4\phantom{-}&0&0&0&0\\
\rep{3{}'\phantom{'}}&3&0&0&3&-1\phantom{-}&1&-1\phantom{-}&-1\phantom{-}\\
\rep{2{}''}&2&-1\phantom{-}&-1\phantom{-}&2&0&0&2&0\\
\rep{2{}'\phantom{'}}&2&1&-1\phantom{-}&-2\phantom{-}&\sigma&0&0
&-\sigma\phantom{-}\\
\rep{1{}'\phantom{'}}&1&1&1&1&-1\phantom{-}&-1\phantom{-}&1&-1\phantom{-}
\end{array}$\qquad$\sigma=\sqrt2$

\subsection*{Icosahedral case
$\III=\bigl\langle\alpha,\beta,\gamma\bigm|\alpha^5=\beta^3=\gamma^2
=\alpha\beta\gamma\bigr\rangle\cong\SL_2(\FF_5)$}
\begin{remark}
$\alpha\mapsto\binom{-1\ {-1}}{\phantom{+}0\ {-1}}$,
$\beta\mapsto\binom{\phantom{+}1\ {\phantom{+}1}}{-1\ {\phantom{+}0}}$,
$\gamma\mapsto\binom{\,0\ {-1}}{\,1 \ {\phantom{+}0}}$ gives
an isomorphism $\III\cong\SL_2(\FF_5)$.
\end{remark}

\setlength{\unitlength}{1cm}
\begin{center}
{\footnotesize
\begin{picture}(7,2)(-5,-1)
\multiput(-5,0)(1,0){8}{\circle{0.2}}
\multiput(-4,0)(2,0){4}{\circle*{0.08}}
\put(0,-1){\circle{0.2}}
\put(-5,0){\makebox(0,0){\tiny$\times$}}
\multiput(-4.9,0)(1,0){7}{\line(1,0){0.8}}
\put(0,-0.9){\line(0,1){0.8}}
\put(-5,0.237){\makebox(0,0)[b]{$\substack{\textstyle\{1\}\\[0.456ex]
\textstyle\rep{1}}$}}
\put(-4,0.237){\makebox(0,0)[b]{$\substack{\textstyle[\alpha]\\[0.456ex]
\textstyle\rep{2}}$}}
\put(-3,0.237){\makebox(0,0)[b]{$\substack{\textstyle[\alpha^2]\\[0.456ex]
\textstyle\rep{3}}$}}
\put(-2,0.237){\makebox(0,0)[b]{$\substack{\textstyle[\alpha^3]\\[0.456ex]
\textstyle\rep{4}}$}}
\put(-1,0.237){\makebox(0,0)[b]{$\substack{\textstyle[\alpha^4]\\[0.456ex]
\textstyle\rep{5}}$}}
\put(0,0.237){\makebox(0,0)[b]{$\substack{\textstyle\{\overline1\}\\[0.456ex]
\textstyle\rep{6}}$}}
\put(1,0.237){\makebox(0,0)[b]{$\substack{\textstyle[\beta^2]\\[0.456ex]
\textstyle\rep{4'}}$}}
\put(2,0.237){\makebox(0,0)[b]{$\substack{\textstyle[\beta]\\[0.456ex]
\textstyle\rep{2'}}$}}
\put(0.237,-1){\makebox(0,0)[l]{$\rep{3'}\ [\gamma]$}}
\end{picture}
}
\end{center}

$\begin{array}[b]{r|ccccccccc}
\multicolumn{1}{p{\cwth}}{}&\multicolumn{1}{p{\cwth}}{}&
\multicolumn{1}{p{\cwth}}{}&\multicolumn{1}{p{\cwth}}{}&
\multicolumn{1}{p{\cwth}}{}&\multicolumn{1}{p{\cwth}}{}&
\multicolumn{1}{p{\cwth}}{}&\multicolumn{1}{p{\cwth}}{}&
\multicolumn{1}{p{\cwth}}{}&\multicolumn{1}{p{\cwth}}{}\\
\substack{\textrm{character}\\\textrm{table}\hfill}
&\{1\}&[\alpha]
&[\alpha^2]&[\alpha^3]
&[\alpha^4]&\{\overline1\}
&[\beta^2]&[\gamma]
&[\beta]\\\hline
\rep{1\phantom{'}}&1&1&1&1&1&1&1&1&1\\
\rep{2\phantom{'}}&2&\tau&-\tau'\phantom{-}&\tau'&-\tau\phantom{-}
&-2\phantom{-}&-1\phantom{-}&0&1\\
\rep{3\phantom{'}}&3&\tau&\tau'&\tau'&\tau&3&0&-1\phantom{-}&0\\
\rep{4\phantom{'}}&4&1&-1\phantom{-}&1&-1\phantom{-}&-4\phantom{-}&1&0
&-1\phantom{-}\\
\rep{5\phantom{'}}&5&0&0&0&0&5&-1\phantom{-}&1&-1\phantom{-}\\
\rep{6\phantom{'}}&6&-1\phantom{-}&1&-1\phantom{-}&1&-6\phantom{-}&0&0&0\\
\rep{4{}'}&4&-1\phantom{-}&-1\phantom{-}&-1\phantom{-}&-1\phantom{-}&4&1&0&1\\
\rep{3{}'}&3&\tau'&\tau&\tau&\tau'&3&0&-1\phantom{-}&0\\
\rep{2{}'}&2&\tau'&-\tau\phantom{-}&\tau&-\tau'\phantom{-}&-2\phantom{-}
&-1\phantom{-}&0&1
\end{array}$\qquad$\begin{array}{@{}r}
\tau=\frac{1+\sqrt5}{2}\\
\tau'=\frac{1-\sqrt5}{2}\end{array}$

\section{More on the three primitive cases}\label{sect}
Let us denote the irreducible representations of $\SU(2)$ (or rather their
classes in the representation ring) by bold letters. So \mbox{\boldmath$i$}
stands for the $i$-dimensional representation
$\mbox{\boldmath$i$}=\Sym^{i-1}(\Rep{2})$ of $\SU(2)$. As before we use sanserif
fonts for representations of $\GrpG$. Addition in the representation ring will
be denoted by $+$ and multiplication usually by $\otimes$. We also write,
e.\,g., $2\cdot\rep5$ for $\rep5+\rep5$. Multiplication in
$\REP(\GrpG)[\![t]\!]$ will be indicated by juxtaposition while we keep the
tensor product sign if no $t$ is in sight.

If we denote the character of $\Rep{2}$ by $\xi+\xi^{-1}$, then
$\mbox{\boldmath$i$}=\Sym^{i-1}(\Rep{2})$ has character
$\sum\limits_{k=0}^{i-1}\xi^{i-1-2k}$ and $\lambda_{-t}\mbox{\boldmath$i$}$
becomes $\prod\limits_{k=0}^{i-1}\bigl(1-\xi^{i-1-2k}t\bigr)$.
Explicitly, the first few $\lambda_{-t}\mbox{\boldmath$i$}$ read as follows.

\begin{align*}
\lambda_{-t}\Rep{1}&=\Rep{1}-\Rep{1}t\\
\lambda_{-t}\Rep{2}&=\Rep{1}-\Rep{2}t+\Rep{1}t^2\\
\lambda_{-t}\Rep{3}&=\Rep{1}-\Rep{3}t+\Rep{3}t^2-\Rep{1}t^3\\
\lambda_{-t}\Rep{4}&=\Rep{1}-\Rep{4}t+(\Rep{1}+\Rep{5})t^2-\Rep{4}t^3+
\Rep{1}t^4\\
\lambda_{-t}\Rep{5}&=\Rep{1}-\Rep{5}t+(\Rep{3}+\Rep{7})t^2-(\Rep{3}+
\Rep{7})t^3+\Rep{5}t^4-\Rep{1}t^5\\
\lambda_{-t}\Rep{6}&=\Rep{1}-\Rep{6}t+(\Rep{1}+\Rep{5}+\Rep{9})t^2
-(\Rep{4}+\Rep{6}+\Rep{10})t^3+(\Rep{1}+\Rep{5}+\Rep{9})t^4-\Rep{6}t^5+
\Rep{1}t^6\\
\lambda_{-t}\Rep{7}&=\Rep{1}-\Rep{7}t+(\Rep{3}+\Rep{7}+\Rep{11})t^2
-(\Rep{1}+\Rep{5}+\Rep{7}+\Rep{9}+\Rep{13})t^3\\
&\phantom{{}=\Rep{1}-\Rep{7}t+(\Rep{3}+\Rep{7}+\Rep{11})t^2}+
(\Rep{1}+\Rep{5}+\Rep{7}+\Rep{9}+\Rep{13})t^4-(\Rep{3}+\Rep{7}+\Rep{11})t^5+
\Rep{7}t^6-\Rep{1}t^7
\end{align*}

Next we shall compute the $\lambda$-polynomials $\lambda_{-t}(\rep i)$ for
each irreducible representation in the binary tetrahedral, octahedral,
and icosahedral cases.

\subsection*{$\TTT$}
$\begin{array}{@{}c||c|c|c|c|c|}
\SU(2)&\Rep{1}&\Rep{2}&\Rep{3}&\Rep{4}&\Rep{5}\\\hline
\TTT&\rep{1}&\rep{2}&\rep{3}&\rep{2'}+\rep{2''}&\rep{3}+\rep{1'}+\rep{1''}
\end{array}$\\[2mm]
Restricting from $\SU(2)$ to $\TTT$ gives the $\lambda$-polynomials
for the representations $\rep{1},\rep{2},\rep{3}$ of~$\TTT$.
From the character table we get $\rep{2'}\otimes\rep{2'}=\rep{3}+\rep{1''}$,
and its $1$-dimensional summand $\bigwedge^2\rep{2'}$ is therefore
$\bigwedge^2\rep{2'}=\rep{1''}$. Similarly we get
$\bigwedge^2\rep{2''}=\rep{1'}$.\\

\noindent
Here is the list of $\lambda$-polynomials for the irreducible representations of
$\TTT$.
\begin{align*}
\lambda_{-t}\rep{1}&=\rep{1}-\rep{1}t\\
\lambda_{-t}\rep{2}&=\rep{1}-\rep{2}t+\rep{1}t^2\\
\lambda_{-t}\rep{3}&=\rep{1}-\rep{3}t+\rep{3}t^2-\rep{1}t^3\\
\lambda_{-t}\rep{2'}&=\rep{1}-\rep{2'}t+\rep{1''}t^2\\
\lambda_{-t}\rep{2''}&=\rep{1}-\rep{2''}t+\rep{1'}t^2\\
\lambda_{-t}\rep{1'}&=\rep{1}-\rep{1'}t\\
\lambda_{-t}\rep{1''}&=\rep{1}-\rep{1''}t
\end{align*}

\subsection*{$\OOO$}
$\begin{array}{@{}c||c|c|c|c|c|c|c|}
\SU(2)&\Rep{1}&\Rep{2}&\Rep{3}&\Rep{4}&\Rep{5}&\Rep{6}&\Rep{7}\\\hline
\OOO&\rep{1}&\rep{2}&\rep{3}&\rep{4}&\rep{3'}+\rep{2''}&\rep{4}+\rep{2'}
&\rep{3}+\rep{3'}+\rep{1'}
\end{array}$\\[2mm]
Restricting from $\SU(2)$ to $\OOO$ gives the $\lambda$-polynomials
for the representations $\rep{1},\dots,\rep{4}$ of~$\OOO$.
From the character table we get $\rep{2'}\otimes\rep{2'}=\rep{3}+\rep{1}$,
so that $\bigwedge^2\rep{2'}=\rep{1}$. From $\rep{2''}\otimes\rep{2''}
=\rep{1}+\rep{2''}+\rep{1'}$, however, one cannot read off
$\bigwedge^2\rep{2''}$. But from
\begin{align*}
\lambda_{-t}(\Rep{7}|_\OOO)&=\lambda_{-t}(\rep{3}+\rep{3'}+\rep{1'})
=\lambda_{-t}\rep{3}\,\lambda_{-t}\rep{3'}\,\lambda_{-t}\rep{1'}\quad
\mbox{and}\\
\lambda_{-t}(\Rep{5}|_\OOO)&=\lambda_{-t}(\rep{3'}+\rep{2''})
=\lambda_{-t}\rep{3'}\,\lambda_{-t}\rep{2''}
\end{align*}
we first get $\bigwedge^3\rep{3'}=\rep{1'}$ and $\bigwedge^2\rep{2''}=\rep{1'}$.
Looking at the $t^2$ coefficient of the second equation above we see that
$$\Rep3|_\OOO+\Rep7|_\OOO=\rep3+\rep3+\rep{3'}+\rep{1'}
=\rep{1'}+\rep{3'}\otimes\rep{2''}+\textstyle\bigwedge^2\rep{3'}$$
and since $\rep{3'}\otimes\rep{2''}=\rep3+\rep{3'}$ we get
$\bigwedge^2\rep{3'}=\rep3$.\\

\noindent
Here is the list of $\lambda$-polynomials for the irreducible representations of
$\OOO$.
\begin{align*}
\lambda_{-t}\rep{1}&=\rep{1}-\rep{1}t\\
\lambda_{-t}\rep{2}&=\rep{1}-\rep{2}t+\rep{1}t^2\\
\lambda_{-t}\rep{3}&=\rep{1}-\rep{3}t+\rep{3}t^2-\rep{1}t^3\\
\lambda_{-t}\rep{4}&=\rep{1}-\rep{4}t+(\rep{1}+\rep{3'}+\rep{2''})t^2-\rep{4}t^3
+\rep{1}t^4\\
\lambda_{-t}\rep{3'}&=\rep{1}-\rep{3'}t+\rep{3}t^2-\rep{1'}t^3\\
\lambda_{-t}\rep{2''}&=\rep{1}-\rep{2''}t+\rep{1'}t^2\\
\lambda_{-t}\rep{2'}&=\rep{1}-\rep{2'}t+\rep{1}t^2\\
\lambda_{-t}\rep{1'}&=\rep{1}-\rep{1'}t
\end{align*}

\subsection*{$\III$}
$\begin{array}{@{}c||c|c|c|c|c|c|c|c|c|c|}
\SU(2)&\Rep{1}&\Rep{2}&\Rep{3}&\Rep{4}&\Rep{5}&\Rep{6}&\Rep{7}
&\Rep{8}&\Rep{9}&\Rep{10}\\\hline
\III&\rep{1}&\rep{2}&\rep{3}&\rep{4}&\rep{5}&\rep{6}&\rep{4'}+\rep{3'}
&\rep{6}+\rep{2'}&\rep{5}+\rep{4'}&\rep{4}+\rep{6}
\end{array}$\\[2mm]
Restricting from $\SU(2)$ to $\III$ gives the $\lambda$-polynomials
for the representations $\rep{1},\dots,\rep{6}$ of~$\III$.
The $\lambda$-polynomials for $\rep{3'}$ and $\rep{2'}$ are obvious, too.
Only the $t^2$ coefficient in $\lambda_{-t}\rep{4'}$ is not totally obvious.
For this, we can use the universal formula
$$\textstyle
\bigwedge^2(x\otimes y)=\bigwedge^2x\otimes y^{\otimes2}+
x^{\otimes2}\otimes\bigwedge^2y-2\cdot\bigwedge^2x\otimes\bigwedge^2y$$
and apply it for $\rep{4'}=\rep{2}\otimes\rep{2'}$.\\

\noindent
Here is the list of $\lambda$-polynomials for the irreducible representations of
$\III$.
\begin{align*}
\lambda_{-t}\rep{1}&=\rep{1}-\rep{1}t\\
\lambda_{-t}\rep{2}&=\rep{1}-\rep{2}t+\rep{1}t^2\\
\lambda_{-t}\rep{3}&=\rep{1}-\rep{3}t+\rep{3}t^2-\rep{1}t^3\\
\lambda_{-t}\rep{4}&=\rep{1}-\rep{4}t+(\rep{1}+\rep{5})t^2-\rep{4}t^3+
\rep{1}t^4\\
\lambda_{-t}\rep{5}&=\rep{1}-\rep{5}t+(\rep{3}+\rep{4'}+\rep{3'})t^2
-(\rep{3}+\rep{4'}+\rep{3'})t^3+\rep{5}t^4-\rep{1}t^5\\
\lambda_{-t}\rep{6}&=\rep{1}-\rep{6}t+(\rep{1}+2\cdot\rep{5}+\rep{4'})t^2
-(2\cdot\rep{4}+2\cdot\rep{6})t^3+(\rep{1}+2\cdot\rep{5}+\rep{4'})t^4
-\rep{6}t^5+\rep{1}t^6\\
\lambda_{-t}\rep{4'}&=\rep{1}-\rep{4'}t+(\rep{3}+\rep{3'})t^2-\rep{4'}t^3+
\rep{1}t^4\\
\lambda_{-t}\rep{3'}&=\rep{1}-\rep{3'}t+\rep{3'}t^2-\rep{1}t^3\\
\lambda_{-t}\rep{2'}&=\rep{1}-\rep{2'}t+\rep{1}t^2
\end{align*}

We are now ready to compute the Poincar\'e series ($\GrpG$ is fixed and
suppressed in the notation)
$$P_{\rep i,\rep j}(t):=
\sum_{n=0}^\infty\dim\Hom_\GrpG\bigl(\rep i,\Sym^n(\rep j)\bigr)\cdot t^n$$
generalizing $P_{\rep i}(t)=P_{\rep i,\rep2}(t)$ from
Definition~\ref{Poincareseries}. To do so one has to invert the polynomial
$\lambda_{-t}(\rep j)$ in the power series ring $\REP(\GrpG)[\![t]\!]$
(see~(\ref{lambdasigma})). This is best done by regarding
$\lambda_{-t}(\rep j)$ as a multiplication operator in
$\operatorname{Aut}_{\ZZ[\![t]\!]}\bigl(\REP(\GrpG)[\![t]\!]\bigr)$ and by
inverting its matrix (with respect to the basis $\Ghat$). The column
corresponding to the trivial representation $\rep1$ of the inverted matrix
is then the vector $\bigl(P_{\rep i,\rep j}(t)\bigr)_{\rep i\in\Ghat}$ we
are looking for.

Note how the Poincar\'e series $P_{\rep i,\rep3}(t)$ is determined by
$P_{\rep i}(t)=P_{\rep i,\rep2}(t)$.
$$P_{\rep i,\rep3}(t)=\begin{cases}
\dfrac{P_{\rep i}(t^{1/2})}{1-t^2}&\mbox{if $\rep i$ is not spinorial,}\\
\ \mathstrut0&\mbox{if $\rep i$ is spinorial}.
\end{cases}$$
Here is one way of proving this. Let $\psi^2$ be the Adams square which
doubles the $\SU(2)$ weights. If $\pm\varpi$ are the weights of $\Rep2$,
then those of $\Rep3$ are $2\varpi,0,-2\varpi$, so that $\Rep3=\psi^2(\Rep2)+
\Rep1$. More generally, note that
$\psi^2\Sym^{n-2i}(\Rep2)=\Sym^{2n-4i}(\Rep2)-\psi^2\Sym^{n-2i-1}(\Rep2)$.
We get
\begin{align*}
\Sym^n(\Rep3)&=\Sym^n\bigl(\psi^2(\Rep2)+\Rep1\bigr)=
\sum_{k=0}^n\Sym^{n-k}\bigl(\psi^2(\Rep2)\bigr)\otimes\Sym^k(\Rep1)=
\sum_{k=0}^n\psi^2\Sym^{n-k}(\Rep2)\\
&=\sum_{i=0}^{\lfloor n/2\rfloor}\bigl(\Sym^{2n-4i}(\Rep2)-
\psi^2\Sym^{n-2i-1}(\Rep2)\bigr)+
\sum_{i=0}^{\lfloor n/2\rfloor}\psi^2\Sym^{n-2i-1}(\Rep2)
=\sum_{i=0}^{\lfloor n/2\rfloor}\Sym^{2n-4i}(\Rep2).
\end{align*}
Now we can restrict from $\SU(2)$ to $\GrpG$ and compute
\begin{align*}
P_{\rep i,\rep3}(t)&=\sum_{n=0}^\infty\dim\Hom_\GrpG\bigl(\rep i,\Sym^n(\rep3)
\bigr)
\cdot t^n=\sum_{n=0}^\infty\sum_{k=0}^{\lfloor n/2\rfloor}
\dim\Hom_\GrpG\bigl(\rep i,\Sym^{2n-4k}(\Rep2)\bigr)
\cdot t^n\\
&=\sum_{k=0}^\infty\sum_{n=2k}^\infty
\dim\Hom_\GrpG\bigl(\rep i,\Sym^{2n-4k}(\rep2)\bigr)
\cdot t^{n-2k}\cdot t^{2k}\\
&=\sum_{k=0}^\infty t^{2k}\sum_{n=0}^\infty
\dim\Hom_\GrpG\bigl(\rep i,\Sym^n(\rep2)\bigr)\cdot t^{n/2}
=\frac{1}{1-t^2}\,P_{\rep i}(t^{1/2})
\end{align*}
assuming that $\rep i$ is not spinorial, so that $\Hom_\GrpG\bigl(\rep i,
\Sym^{\operatorname{odd}}(\Rep2)\bigr)=0$.

It would be easy to list all the Poincar\'e series $P_{\rep i,\rep j}(t)$.
Let it suffice to show the Poincar\'e series $P_{\rep1,\rep j}(t)$
for the invariant rings $\Sym^\ast(\rep j)^\GrpG$. 

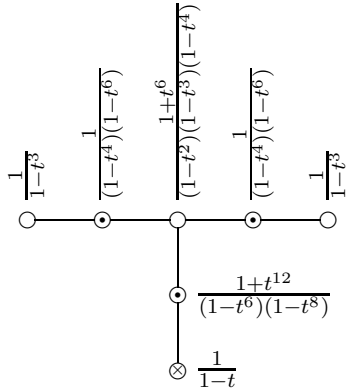
\setlength{\unitlength}{1cm}
\subsection*{$\GrpG=\TTT$}
\begin{center}
{\small
\begin{picture}(4,6)(-2,-2)
\multiput(-2,0)(1,0){5}{\circle{0.2}}
\multiput(-1,0)(2,0){2}{\circle*{0.08}}
\multiput(0,-2)(0,1){2}{\circle{0.2}}
\put(0,-1){\circle*{0.08}}
\put(0,-2){\makebox(0,0){\tiny$\times$}}
\multiput(-1.9,0)(1,0){4}{\line(1,0){0.8}}
\multiput(0,-1.9)(0,1){2}{\line(0,1){0.8}}
\put(-2,0.237){\makebox(0,0)[b]{\begin{turn}{90}$\frac{1}{1-t^3}$\end{turn}}}
\put(-1,0.237){\makebox(0,0)[b]{\begin{turn}{90}$\frac{1}{(1-t^4)(1-t^6)}$
\end{turn}}}
\put(0,0.237){\makebox(0,0)[b]{\begin{turn}{90}$\frac{1+t^6}{(1-t^2)(1-t^3)
(1-t^4)}$\end{turn}}}
\put(1,0.237){\makebox(0,0)[b]{\begin{turn}{90}$\frac{1}{(1-t^4)(1-t^6)}$
\end{turn}}}
\put(2,0.237){\makebox(0,0)[b]{\begin{turn}{90}$\frac{1}{1-t^3}$\end{turn}}}
\put(0.237,-1){\makebox(0,0)[l]{$\frac{1+t^{12}}{(1-t^6)(1-t^8)}$}}
\put(0.237,-2){\makebox(0,0)[l]{$\frac{1}{1-t}$}}
\end{picture}
}
\end{center}

\setlength{\unitlength}{1cm}
\subsection*{$\GrpG=\OOO$}
\begin{center}
{\small
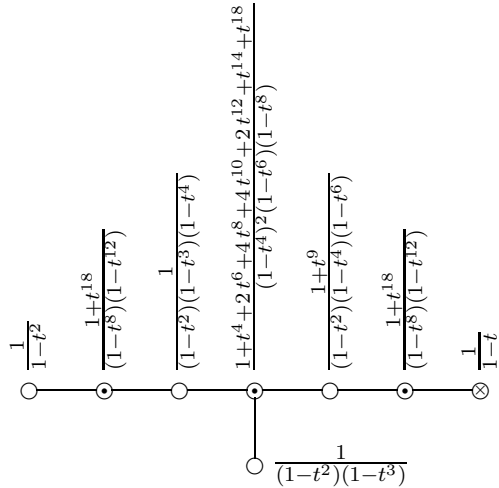
\begin{picture}(6,7)(-3,-1)
\multiput(-3,0)(1,0){7}{\circle{0.2}}
\multiput(-2,0)(2,0){3}{\circle*{0.08}}
\put(0,-1){\circle{0.2}}
\put(3,0){\makebox(0,0){\tiny$\times$}}
\multiput(-2.9,0)(1,0){6}{\line(1,0){0.8}}
\put(0,-0.9){\line(0,1){0.8}}
\put(-3,0.237){\makebox(0,0)[b]{\begin{turn}{90}$\frac{1}{1-t^2}$\end{turn}}}
\put(-2,0.237){\makebox(0,0)[b]{\begin{turn}{90}$\frac{1+t^{18}}{(1-t^8)
(1-t^{12})}$\end{turn}}}
\put(-1,0.237){\makebox(0,0)[b]{\begin{turn}{90}$\frac{1}{(1-t^2)(1-t^3)
(1-t^4)}$\end{turn}}}
\put(0,0.237){\makebox(0,0)[b]{\begin{turn}{90}$\frac{
1+t^4+2\,t^6+4\,t^8+4\,t^{10}+2\,t^{12}+t^{14}+t^{18}}{(1-t^4)^2(1-t^6)
(1-t^8)}$\end{turn}}}
\put(1,0.237){\makebox(0,0)[b]{\begin{turn}{90}$\frac{1+t^9}{(1-t^2)(1-t^4)
(1-t^6)}$\end{turn}}}
\put(2,0.237){\makebox(0,0)[b]{\begin{turn}{90}$\frac{1+t^{18}}{(1-t^8)
(1-t^{12})}$\end{turn}}}
\put(3,0.237){\makebox(0,0)[b]{\begin{turn}{90}$\frac{1}{1-t}$\end{turn}}}
\put(0.237,-1){\makebox(0,0)[l]{$\frac{1}{(1-t^2)(1-t^3)}$}}
\end{picture}
}
\end{center}

\setlength{\unitlength}{1cm}
\subsection*{$\GrpG=\III$}
\begin{center}
{\small
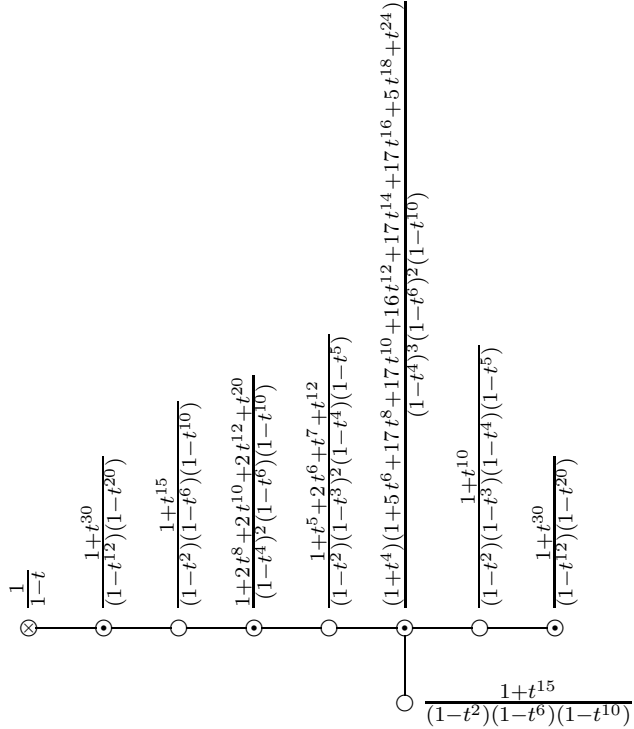
\begin{picture}(7,10)(-5,-1)
\multiput(-5,0)(1,0){8}{\circle{0.2}}
\multiput(-4,0)(2,0){4}{\circle*{0.08}}
\put(0,-1){\circle{0.2}}
\put(-5,0){\makebox(0,0){\tiny$\times$}}
\multiput(-4.9,0)(1,0){7}{\line(1,0){0.8}}
\put(0,-0.9){\line(0,1){0.8}}
\put(-5,0.237){\makebox(0,0)[b]{\begin{turn}{90}$\frac{1}{1-t}$\end{turn}}}
\put(-4,0.237){\makebox(0,0)[b]{\begin{turn}{90}$\frac{1+t^{30}}{(1-t^{12})
(1-t^{20})}$\end{turn}}}
\put(-3,0.237){\makebox(0,0)[b]{\begin{turn}{90}$\frac{1+t^{15}}{(1-t^2)
(1-t^6)(1-t^{10})}$\end{turn}}}
\put(-2,0.237){\makebox(0,0)[b]{\begin{turn}{90}$\frac{1+2\,t^8+2\,t^{10}+
2\,t^{12}+t^{20}}{
(1-t^4)^2(1-t^6)(1-t^{10})}$\end{turn}}}
\put(-1,0.237){\makebox(0,0)[b]{\begin{turn}{90}$\frac{1+t^5+2\,t^6+t^7+
t^{12}}{
(1-t^2)(1-t^3)^2(1-t^4)(1-t^5)}$\end{turn}}}
\put(0,0.237){\makebox(0,0)[b]{\begin{turn}{90}$\frac{(1+t^4)
(1+5\,t^6+17\,t^8+17\,t^{10}+16\,t^{12}+17\,t^{14}+17\,t^{16}+5\,t^{18}+
t^{24})}{
(1-t^4)^3(1-t^6)^2(1-t^{10})}$\end{turn}}}
\put(1,0.237){\makebox(0,0)[b]{\begin{turn}{90}$\frac{1+t^{10}}{(1-t^2)
(1-t^3)(1-t^4)(1-t^5)}$\end{turn}}}
\put(2,0.237){\makebox(0,0)[b]{\begin{turn}{90}$\frac{1+t^{30}}{(1-t^{12})
(1-t^{20})}$\end{turn}}}
\put(0.237,-1){\makebox(0,0)[l]{$\frac{1+t^{15}}{(1-t^2)(1-t^6)(1-t^{10})}$}}
\end{picture}
}
\end{center}

\begin{remark}
Note that in this last case $\GrpG=\III$ one has
$\dim\rep j=-\deg P_{\rep1,\rep j}(t)$ for all $\rep j\in\Ghat$ where the degree
of a rational function is defined to be the degree of the numerator polynomial
minus the degree of the denominator polynomial (and this difference is of course
well-defined).
This example suggests evident classification questions.
\end{remark}

\section{Homomorphisms into finite complex reflection\\groups}\label{hifcrg}
The well-known Chevalley-Shephard-Todd Theorem characterizes the finite complex
reflection groups as those finite groups that have a polynomial algebra as
invariant ring of some faithful complex representation. More precisely, one has
the following theorem.

\begin{theorem}
Let $V$ be a finite-dimensional $K$-vector space and $G\hookrightarrow\GL(V)$ a
faithful representation of a finite group $G$. Assume that $\operatorname{char}K
\nmid|G|$. Then the following are equivalent.
\begin{itemize}
\item[(i)] $G$ is generated by (pseudo-)reflections.
\item[(ii)] $K[V]^G$ is a polynomial algebra.
\end{itemize}
\end{theorem}

The original verification of this theorem for $K=\CC$ by Shephard and Todd
\cite{ST} depended on their classification of finite irreducible complex
reflection groups. A case free (but computational) proof was furnished by
Chevalley \cite{Ch}. A noncomputational proof of this result can be found in
\cite{Sm}.

The aim of this section is to exhibit explicitly, for each of the three
primitive binary polyhedral groups $\GrpG$ and each of their irreducible
representations $\rep i$, a homomorphism from $\GrpG$ to a finite complex
reflection group $G$ such that the reflection representation of $G$ restricts to
the representation $\rep i$. A hint for guessing suitable target groups $G$ for
such homomorphisms $\GrpG\to G$ comes from looking at the exponents in the
denominator of the Poincar\'e series of the invariant rings
$\Sym^\ast(\rep i)^\GrpG$ and comparing them with the degrees for the finite
complex reflection groups as tabulated for instance in \cite{BMR} where one also
finds presentations ``\`a la Coxeter'' for all these groups.

Our results are compiled below. For each primitive binary polyhedral group
$\GrpG$ and each of its irreducible representations $\rep i$ we give
$\GrpG\twoheadrightarrow\overline{\GrpG}\hookrightarrow G$ for a complex
reflection group $G$ in the Shephard-Todd list. We make this very explicit by
writing down a presentation for $G$ and by giving the matrices for the
reflection representation (unless $G$ is a symmetric group; if $G$ is a
symmetric group then the usual representation by permutation matrices is the
direct sum of the reflection representation and a $1$-dimensional trivial
representation). Finally, the homomorphisms $\GrpG\to G$ are defined by giving
the images of the generators $\alpha$ and $\beta$ of $\GrpG$. One can then
easily check that this gives the correct character values. The computer algebra
system GAP \cite{GAP} was used for some computations.

\newlength{\frepbox}
\settowidth{\frepbox}{$\rep2''$}

\subsection*{Tetrahedral case
$\TTT=\bigl\langle\alpha,\beta,\gamma\bigm|\alpha^3=\beta^3=\gamma^2
=\alpha\beta\gamma\bigr\rangle$}
$\TTT$ has four normal subgroups. The quotient groups look as follows.
\begin{itemize}
\item $\bigl\langle\alpha,\beta,\gamma\bigm|\alpha^3=\beta^3=\gamma^2
=\alpha\beta\gamma\bigr\rangle=\TTT$
\item $\bigl\langle\alpha,\beta,\gamma\bigm|\alpha^3=\beta^3=\gamma^2
=\alpha\beta\gamma=1\bigr\rangle\cong\Al_4$
\item $\bigl\langle\alpha,\beta,\gamma\bigm|\alpha^3=\beta^3=\gamma
=\alpha\beta=1\bigr\rangle\cong\Al_3$
\item $\bigl\langle\alpha,\beta,\gamma\bigm|\alpha=\beta=\gamma=1
\bigr\rangle=1$
\end{itemize}

\noindent
\fbox{\makebox[\frepbox][c]{$\mathstrut\rep1$}}\quad $\TTT\twoheadrightarrow1
\hookrightarrow\Sy_2\cong W(\sfA_1)=\bigl\langle r\bigm|r^2=1\bigr\rangle$
$$\alpha\mapsto1\quad\beta\mapsto1$$

\noindent
\fbox{\makebox[\frepbox][c]{$\mathstrut\rep2$}}\quad $\TTT\hookrightarrow G_{12}
=\bigl\langle r,s,t\bigm|r^2=s^2=t^2=1,\ rstr=strs=trst\bigr\rangle$
$$r\mapsto\Bigl(\begin{smallmatrix}
\phantom{+}0&\varepsilon\\-\varepsilon^3&0\end{smallmatrix}\Bigr)\quad
s\mapsto\tfrac{1}{\sqrt2}\Bigl(\begin{smallmatrix}
\phantom{+}1&\phantom{+}i\\-i&-1\end{smallmatrix}\Bigr)\quad
t\mapsto\tfrac{1}{\sqrt2}\Bigl(\begin{smallmatrix}
1&\phantom{+}1\\1&-1\end{smallmatrix}\Bigr)\qquad(\varepsilon=\exp(\pi i/4))$$
$$\alpha\mapsto rs\quad\beta\mapsto ts$$

\noindent
\fbox{\makebox[\frepbox][c]{$\mathstrut\rep3$}}\quad $\TTT\twoheadrightarrow
\Al_4\hookrightarrow\Sy_4\cong W(\sfA_3)=\bigl\langle r,s,t\bigm|
r^2=s^2=t^2=(rs)^3=(st)^3=(rt)^2=1\bigr\rangle$\\[-1mm]
$$\alpha\mapsto rs\quad\beta\mapsto st$$

\noindent
\fbox{\makebox[\frepbox][c]{$\mathstrut\rep{2'}$}}\quad $\TTT
\stackrel{\cong}{\to}G_4=\bigl\langle r,s\bigm|r^3=s^3=1,\ rsr=srs\bigr\rangle$
$$r\mapsto\Bigl(\begin{smallmatrix}
1&0\\0&\rho\end{smallmatrix}\Bigr)\quad
s\mapsto\tfrac{i}{\sqrt3}\Bigl(\begin{smallmatrix}
\phantom{+}1&-\sqrt2\,\rho^2\\-\sqrt2\,\rho^2&-\rho\end{smallmatrix}\Bigr)
\qquad(\rho=\exp(2\pi i/3))$$
$$\alpha\mapsto r^2s^2\quad\beta\mapsto rs$$

\noindent
\fbox{\makebox[\frepbox][c]{$\mathstrut\rep{2''}$}}\quad $\TTT
\stackrel{\cong}{\to}G_4=\bigl\langle r,s\bigm|r^3=s^3=1,\ rsr=srs\bigr\rangle$
$$\alpha\mapsto rs\quad\beta\mapsto r^2s^2$$

\noindent
\fbox{\makebox[\frepbox][c]{$\mathstrut\rep{1'}$}}\quad $\TTT\twoheadrightarrow
\Al_3\cong G_3(3)=\bigl\langle r\bigm|r^3=1\bigr\rangle$, $r\mapsto\rho=
\exp(2\pi i/3)$
$$\alpha\mapsto r^2\quad\beta\mapsto r$$

\noindent
\fbox{\makebox[\frepbox][c]{$\mathstrut\rep{1''}$}}\quad $\TTT\twoheadrightarrow
\Al_3\cong G_3(3)=\bigl\langle r\bigm|r^3=1\bigr\rangle$
$$\alpha\mapsto r\quad\beta\mapsto r^2$$

\subsection*{Octahedral case
$\OOO=\bigl\langle\alpha,\beta,\gamma\bigm|\alpha^4=\beta^3=\gamma^2
=\alpha\beta\gamma\bigr\rangle$}
$\OOO$ has five normal subgroups. The quotient groups look as follows.
\begin{itemize}
\item $\bigl\langle\alpha,\beta,\gamma\bigm|\alpha^4=\beta^3=\gamma^2
=\alpha\beta\gamma\bigr\rangle=\OOO$
\item $\bigl\langle\alpha,\beta,\gamma\bigm|\alpha^4=\beta^3=\gamma^2
=\alpha\beta\gamma=1\bigr\rangle\cong\Sy_4$
\item $\bigl\langle\alpha,\beta,\gamma\bigm|\alpha^2=\beta^3=\gamma^2
=\alpha\beta\gamma=1\bigr\rangle\cong\Sy_3$
\item $\bigl\langle\alpha,\beta,\gamma\bigm|\alpha^2=\beta=\gamma^2
=\alpha\gamma=1\bigr\rangle\cong\Sy_2$
\item $\bigl\langle\alpha,\beta,\gamma\bigm|\alpha=\beta=\gamma=1
\bigr\rangle=1$
\end{itemize}

\noindent
\fbox{\makebox[\frepbox][c]{$\mathstrut\rep1$}}\quad $\OOO\twoheadrightarrow1
\hookrightarrow\Sy_2\cong W(\sfA_1)=\bigl\langle r\bigm|r^2=1\bigr\rangle$
$$\alpha\mapsto1\quad\beta\mapsto1$$

\noindent
\fbox{\makebox[\frepbox][c]{$\mathstrut\rep2$}}\quad $\OOO\hookrightarrow
G_{13}=\bigl\langle r,s,t\bigm|
r^2=s^2=t^2=1,\ rstrs=trstr,\ strs=trst\bigr\rangle$
$$r\mapsto\Bigl(\begin{smallmatrix}
-1&0\\\phantom{+}0&1\end{smallmatrix}\Bigr)\quad
s\mapsto\tfrac{1}{\sqrt2}\Bigl(\begin{smallmatrix}
-1&-i\\\phantom{+}i&\phantom{+}1\end{smallmatrix}\Bigr)\quad
t\mapsto\tfrac{1}{\sqrt2}\Bigl(\begin{smallmatrix}
-1&1\\\phantom{+}1&1\end{smallmatrix}\Bigr)$$
$$\alpha\mapsto rs\quad\beta\mapsto ts$$

\noindent
\fbox{\makebox[\frepbox][c]{$\mathstrut\rep3$}}\quad $\OOO\twoheadrightarrow
\Sy_4\hookrightarrow W(\sfB_3)=\bigl\langle r,s,t\bigm|
r^2=s^2=t^2=(rs)^4=(st)^3=(rt)^2=1\bigr\rangle$
$$r\mapsto\Bigl(\begin{smallmatrix}
-1&0&0\\\phantom{+}0&1&0\\\phantom{+}0&0&1
\end{smallmatrix}\Bigr)\quad
s\mapsto\Bigl(\begin{smallmatrix}
0&1&0\\1&0&0\\0&0&1\end{smallmatrix}\Bigr)\quad
t\mapsto\Bigl(\begin{smallmatrix}
1&0&0\\0&0&1\\0&1&0\end{smallmatrix}\Bigr)$$
$$\alpha\mapsto rs\quad\beta\mapsto st$$

\noindent
\fbox{\makebox[\frepbox][c]{$\mathstrut\rep4$}}\quad $\OOO\hookrightarrow
G(4,2,4)=\left\langle r,s,t,u,v\Biggm|
\begin{array}{@{\,}l@{\,}}
r^2=s^2=t^2=u^2=v^2=1,\,rst=str=trs\\
(ru)^2=(rv)^2=1\\
(su)^3=(tu)^3=(uv)^3=(sv)^2=(tv)^2=1
\end{array}\right\rangle$
$$r\mapsto\left(\begin{smallmatrix}
-1&0&0&0\\\phantom{+}0&1&0&0\\\phantom{+}0&0&1&0\\\phantom{+}0&0&0&1
\end{smallmatrix}\right)\quad
s\mapsto\left(\begin{smallmatrix}
0&-i&\phantom{+}0&\phantom{+}0\\
i&\phantom{+}0&\phantom{+}0&\phantom{+}0\\
0&\phantom{+}0&\phantom{+}1&\phantom{+}0\\
0&\phantom{+}0&\phantom{+}0&\phantom{+}1
\end{smallmatrix}\right)\quad
t\mapsto\left(\begin{smallmatrix}
0&1&0&0\\1&0&0&0\\0&0&1&0\\0&0&0&1\end{smallmatrix}\right)\quad
u\mapsto\left(\begin{smallmatrix}
1&0&0&0\\0&0&1&0\\0&1&0&0\\0&0&0&1\\\end{smallmatrix}\right)\quad
v\mapsto\left(\begin{smallmatrix}
1&0&0&0\\0&1&0&0\\0&0&0&1\\0&0&1&0\end{smallmatrix}\right)$$
$$\alpha\mapsto tusrsv\quad\beta\mapsto(rstuv)^2$$

\noindent
\fbox{\makebox[\frepbox][c]{$\mathstrut\rep{3'}$}}\quad $\OOO\twoheadrightarrow
\Sy_4\cong W(\sfA_3)=\bigl\langle r,s,t\bigm|
r^2=s^2=t^2=(rs)^3=(st)^3=(rt)^2=1\bigr\rangle$
$$\alpha\mapsto rst\quad\beta\mapsto ts$$

\noindent
\fbox{\makebox[\frepbox][c]{$\mathstrut\rep{2''}$}}\quad $\OOO\twoheadrightarrow
\Sy_3\cong W(\sfA_2)=\bigl\langle r,s\bigm|r^2=s^2=(rs)^3=1\bigr\rangle$
$$\alpha\mapsto r\quad\beta\mapsto rs$$

\noindent
\fbox{\makebox[\frepbox][c]{$\mathstrut\rep{2'}$}}\quad $\OOO\hookrightarrow
G_{13}=\bigl\langle r,s,t\bigm|
r^2=s^2=t^2=1,\ rstrs=trstr,\ strs=trst\bigr\rangle$
$$\alpha\mapsto(rs)^3\quad\beta\mapsto st$$

\noindent
\fbox{\makebox[\frepbox][c]{$\mathstrut\rep{1'}$}}\quad $\OOO\twoheadrightarrow
\Sy_2\cong W(\sfA_1)=\bigl\langle r\bigm|r^2=1\bigr\rangle$
$$\alpha\mapsto r\quad\beta\mapsto1$$

\subsection*{Icosahedral case
$\III=\bigl\langle\alpha,\beta,\gamma\bigm|\alpha^5=\beta^3=\gamma^2
=\alpha\beta\gamma\bigr\rangle$}
$\III$ has three normal subgroups. The quotient groups look as follows.
\begin{itemize}
\item $\bigl\langle\alpha,\beta,\gamma\bigm|\alpha^5=\beta^3=\gamma^2
=\alpha\beta\gamma\bigr\rangle=\III$
\item $\bigl\langle\alpha,\beta,\gamma\bigm|\alpha^5=\beta^3=\gamma^2
=\alpha\beta\gamma=1\bigr\rangle\cong\Al_5$
\item $\bigl\langle\alpha,\beta,\gamma\bigm|\alpha=\beta=\gamma=1
\bigr\rangle=1$
\end{itemize}

\noindent
\fbox{\makebox[\frepbox][c]{$\mathstrut\rep1$}}\quad $\III\twoheadrightarrow1
\hookrightarrow\Sy_2\cong W(\sfA_1)=\bigl\langle r\bigm|r^2=1\bigr\rangle$
$$\alpha\mapsto1\quad\beta\mapsto1$$

\noindent
\fbox{\makebox[\frepbox][c]{$\mathstrut\rep2$}}\quad $\III\hookrightarrow
G_{22}=\bigl\langle r,s,t\bigm|
r^2=s^2=t^2=1,\,rstrsr=(trs)^2,\,strs=trst\bigr\rangle$
$$r\mapsto\Bigl(\begin{smallmatrix}
0&-i\\i&\phantom{+}0\end{smallmatrix}\Bigr)\quad
s\mapsto\Bigl(\begin{smallmatrix}
\phantom{+}0&\phantom{+}i\eta^3\\-i\eta^2&\phantom{+}0\end{smallmatrix}\Bigr)
\quad
t\mapsto\tfrac{i}{\sqrt5}\Bigl(\begin{smallmatrix}\eta^2-\eta^3&-1+\eta^2\\
1-\eta^3&-\eta^2+\eta^3\end{smallmatrix}\Bigr)\qquad
(\eta=\exp(2\pi i/5))$$
$$\alpha\mapsto rs\quad\beta\mapsto ts$$

\noindent
\fbox{\makebox[\frepbox][c]{$\mathstrut\rep3$}}\quad $\III\twoheadrightarrow
\Al_5\hookrightarrow G_{23}=W(\mathsf H_3)=\bigl\langle r,s,t\bigm|
r^2=s^2=t^2=(rs)^5=(st)^3=(rt)^2=1\bigr\rangle$
$$r\mapsto\tfrac12\Bigl(\begin{smallmatrix}
\tau'&\tau&1\\\tau&1&\tau'\\1&\tau'&\tau
\end{smallmatrix}\Bigr)\quad
s\mapsto\Bigl(\begin{smallmatrix}
-1&0&0\\\phantom{+}0&1&0\\\phantom{+}0&0&1
\end{smallmatrix}\Bigr)\quad
t\mapsto\tfrac12\Bigl(\begin{smallmatrix}
\phantom{+}1&\phantom{+}\tau'&-\tau\\
\phantom{+}\tau'&\phantom{+}\tau&-1\\
-\tau&-1&\phantom{+}\tau'
\end{smallmatrix}\Bigr)
\quad\bigl(\tau=\tfrac{1+\sqrt5}{2}\,,\,\tau'=1-\tau\bigr)$$
$$\alpha\mapsto rs\quad\beta\mapsto st$$

\noindent
\fbox{\makebox[\frepbox][c]{$\mathstrut\rep4$}}\quad $\III\hookrightarrow
G_{29}=\biggl\langle r,s,t,u\biggm|
\begin{array}{@{\,}l@{\,}}
r^2=s^2=t^2=u^2=(rs)^3=(rt)^2=(ru)^2=1\\
(st)^4=(su)^3=(tu)^3=1,\,(ust)^2=(stu)^2
\end{array}\biggr\rangle$
$$r\mapsto\left(\begin{smallmatrix}
-1&0&0&0\\\phantom{+}0&1&0&0\\\phantom{+}0&0&1&0\\\phantom{+}0&0&0&1
\end{smallmatrix}\right)\quad
s\mapsto\tfrac12\left(\begin{smallmatrix}
\phantom{+}1&-1&\phantom{+}i&\phantom{+}i\\
-1&\phantom{+}1&\phantom{+}i&\phantom{+}i\\
-i&-i&\phantom{+}1&-1\\
-i&-i&-1&\phantom{+}1
\end{smallmatrix}\right)\quad
t\mapsto\left(\begin{smallmatrix}
1&\phantom{+}0&\phantom{+}0&\phantom{+}0\\
0&\phantom{+}1&\phantom{+}0&\phantom{+}0\\
0&\phantom{+}0&\phantom{+}0&-1\\
0&\phantom{+}0&-1&\phantom{+}0
\end{smallmatrix}\right)\quad
u\mapsto\left(\begin{smallmatrix}
1&0&0&0\\0&0&1&0\\0&1&0&0\\0&0&0&1
\end{smallmatrix}\right)$$
$$\alpha\mapsto(rstu)^2\quad\beta\mapsto rstsrustst$$

\noindent
\fbox{\makebox[\frepbox][c]{$\mathstrut\rep5$}}\quad $\III\twoheadrightarrow
\Al_5\hookrightarrow\Sy_6\cong W(\sfA_5)$\\
$\phantom{\fbox{$\rep5$}\quad\III\twoheadrightarrow\Al_5\hookrightarrow}{}=
\left\langle r,s,t,u,v\Biggm|
\begin{array}{@{\,}l@{\,}}
r^2=s^2=t^2=u^2=v^2=1\\
(rs)^3=(st)^3=(tu)^3=(uv)^3=1\\
(rt)^2=(ru)^2=(rv)^2=(su)^2=(sv)^2=(tv)^2=1
\end{array}\right\rangle$
$$\alpha\mapsto rstu\quad\beta\mapsto sruv$$
This may be used to recover the classical fact that the outer automorphism group
of the symmetric group of degree $6$ is nontrivial.\\

\noindent
\fbox{\makebox[\frepbox][c]{$\mathstrut\rep6$}}\quad {\LARGE$\mbox{\normalsize
$\III\hookrightarrow G(4,4,6)$}$\\
$\mbox{\normalsize${\phantom{\fbox{$\rep6$}\quad\III\hookrightarrow}}{}=$}
\left\langle
\mbox{\normalsize$\!\!r,s,t,u,v,w\!\!$}\Biggm|\mbox{\normalsize$
\begin{array}{@{\!}l@{\,}}
r^2=s^2=t^2=u^2=v^2=w^2=1\\
(rs)^4=(rt)^3=(ru)^2=(rv)^2=(rw)^2=1\\
(st)^3=1,\,(trs)^2=(rst)^2\\
(su)^2=(sv)^2=(sw)^2=(tv)^2=(tw)^2=(uw)^2=1\\
(tu)^3=(uv)^3=(vw)^3=1
\end{array}$}
\right\rangle$}
$$r\mapsto\left(\begin{smallmatrix}
\phantom{+}0&i&0&0&0&0\\
-i&0&0&0&0&0\\
\phantom{+}0&0&1&0&0&0\\
\phantom{+}0&0&0&1&0&0\\
\phantom{+}0&0&0&0&1&0\\
\phantom{+}0&0&0&0&0&1
\end{smallmatrix}\right)\quad
s\mapsto\left(\begin{smallmatrix}
0&1&0&0&0&0\\1&0&0&0&0&0\\
0&0&1&0&0&0\\0&0&0&1&0&0\\
0&0&0&0&1&0\\0&0&0&0&0&1
\end{smallmatrix}\right)\quad
t\mapsto\left(\begin{smallmatrix}
1&0&0&0&0&0\\0&0&1&0&0&0\\
0&1&0&0&0&0\\0&0&0&1&0&0\\
0&0&0&0&1&0\\0&0&0&0&0&1
\end{smallmatrix}\right)\quad
\cdots\quad
w\mapsto\left(\begin{smallmatrix}
1&0&0&0&0&0\\0&1&0&0&0&0\\
0&0&1&0&0&0\\0&0&0&1&0&0\\
0&0&0&0&0&1\\0&0&0&0&1&0
\end{smallmatrix}\right)$$
$$[\mbox{$u,v$ are mapped to the permutation matrices for the
transpositions $(3\,4),(4\,5)$}]$$
$$\alpha\mapsto(rstuvw)^2\quad
\beta\mapsto srtsrutsrtuvutrwvu$$

\noindent
\fbox{\makebox[\frepbox][c]{$\mathstrut\rep{4'}$}}\quad $\III\twoheadrightarrow\Al_5\hookrightarrow
\Sy_5\cong W(\sfA_4)$\\
$\phantom{\fbox{$\rep{4'}$}\quad\III\twoheadrightarrow\Al_5\hookrightarrow}{}=\biggl\langle r,s,t,u
\biggm|
\begin{array}{@{\,}l@{\,}}
r^2=s^2=t^2=u^2=(rs)^3=(st)^3=(tu)^3=1\\
(rt)^2=(ru)^2=(su)^2=1
\end{array}\biggr\rangle$
$$\alpha\mapsto rstu\quad\beta\mapsto stsr$$

\noindent
\fbox{\makebox[\frepbox][c]{$\mathstrut\rep{3'}$}}\quad $\III\twoheadrightarrow
\Al_5\hookrightarrow G_{23}=W(\mathsf H_3)=\bigl\langle r,s,t\bigm|
r^2=s^2=t^2=(rs)^5=(st)^3=(rt)^2=1\bigr\rangle$
$$\alpha\mapsto (rs)^3\quad\beta\mapsto rtsrsrts$$

\noindent
\fbox{\makebox[\frepbox][c]{$\mathstrut\rep{2'}$}}\quad $\III\hookrightarrow
G_{22}=\bigl\langle r,s,t\bigm|
r^2=s^2=t^2=1,\,rstrsr=(trs)^2,\,strs=trst\bigr\rangle$
$$\alpha\mapsto (rs)^3\quad\beta\mapsto(sr)^2(tr)^2$$

\end{document}